\documentclass{article}
\usepackage[utf8]{inputenc}
\usepackage{vmargin}
\usepackage{amsthm}
\usepackage{amssymb}
\usepackage{amsmath}
\usepackage{stmaryrd}
\usepackage[dvipsnames ,x11names ,svgnames]{xcolor}
\usepackage[none]{hyphenat}
\usepackage{biblatex}
\usepackage{parskip}
\usepackage{graphicx}
\usepackage{listings}
\usepackage{marvosym}
\usepackage{mathrsfs}
\setpapersize{A4}
\setmargins{2 cm}       
{1.5cm}                 
{16.5cm}                
{23.42cm}               
{10pt}                  
{1cm}                   
{0pt}                   
{2cm}                   

\title{Crochet Representation of the Lobachevskian Surface} 
\author{Isabella Estrada Reyes$^\dag$ and Adriana Mejia Casta\~no$^{\dag,1}$ }
\date{}
\begin{document}
\maketitle

\stepcounter{footnote}
\footnotetext{{Correspoding author}}

\begin{abstract}
Beginning the study of non-Euclidean geometries, physical models or representations, such as crochet ones, provide a tangible portrayal of these advanced mathematical concepts. However, their connection to local Euclidean surfaces still needs further investigation. This work aims to explore how the characteristics of crochet models relate to non-Euclidean concepts by providing a parametrization of such surfaces.\\[0.2cm]
\textsl{MSC}: 51L20, 51M10\\
\textsl{Keywords}: Gaussian curvature, hyperbolic geometry, Lobachevskian surface, locally Euclidean surface, Mathematical visualization.
\end{abstract}

\section{Introduction}
The Lobachevskian surface emerges as a field of study that deviates from Euclidean geometry, challenging our intuition and traditional geometric understanding. Fundamental concepts in this surface are often characterized by their non-intuitive properties, owing to the constant negative curvature, as posited in \cite{C}. 
Physical models, known as \emph{representation systems}, have proven to be effective resources enabling a tangible and manipulable representation of abstract ideas \cite{CO}. Crochet models (Figure \ref{coral}), in particular, have garnered attention due to their ability to transform mathematical elements into attractive and concrete three-dimensional objects. However, it remains to be investigated to what extent and how these crochet models can contribute to understanding the mathematical concepts inherent in non-Euclidean geometries.\\

\begin{figure}[h]
        \centering
        {\includegraphics[width=0.3\textwidth]{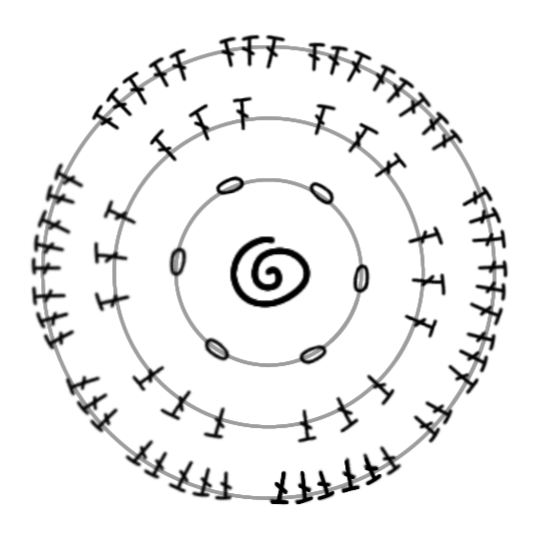}}
        \hspace{2cm}
        {\includegraphics[width=0.3\textwidth]{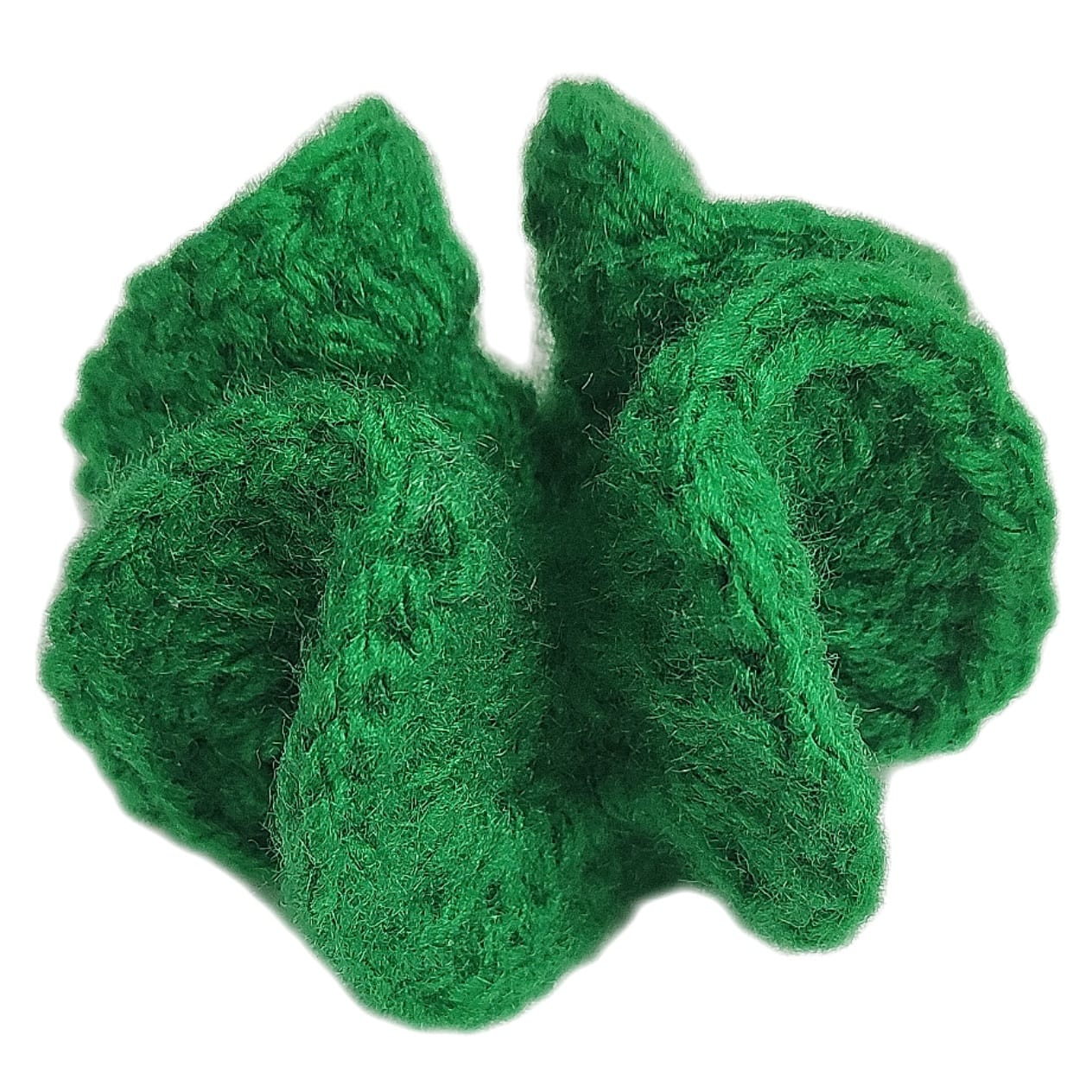}}
        \caption{Crochet diagram (left) and knitted coral (right).}
        \label{coral}
\end{figure}

This work aims to address the following issue: In what ways are the geometric characteristics of crochet models connected to the fundamental mathematical concepts of the Lobachevskian surface, and how can this relationship enhance understanding and intuition? 
\medbreak

Acording to \cite[Section 15]{G}, the \emph{Lobachevskian surface} is the surface of one sheet of the two-sheeted hyperboloid in the pseudo-Euclidean space.
Various models within the literature facilitate the understanding and visualization of this surface. Among these, the most prominent are the Poincar\'e disk and semiplane, both of which are 2-dimensional models.The Poincar\'e disk is defined as the conformal projection of the Lobachevskian surface onto the unitary circle. In this representation, the circumference symbolizes the infinite boundary of the Lobachevskian surface, with the interior corresponding to the Lobachevskian surface itself. Alternatively, an inversion transforms the Poincar\'e disk into the Poincar\'e semi-plane, also preserving angles.

\begin{figure}[h]
    \centering
    {\includegraphics[width=0.2\textwidth]{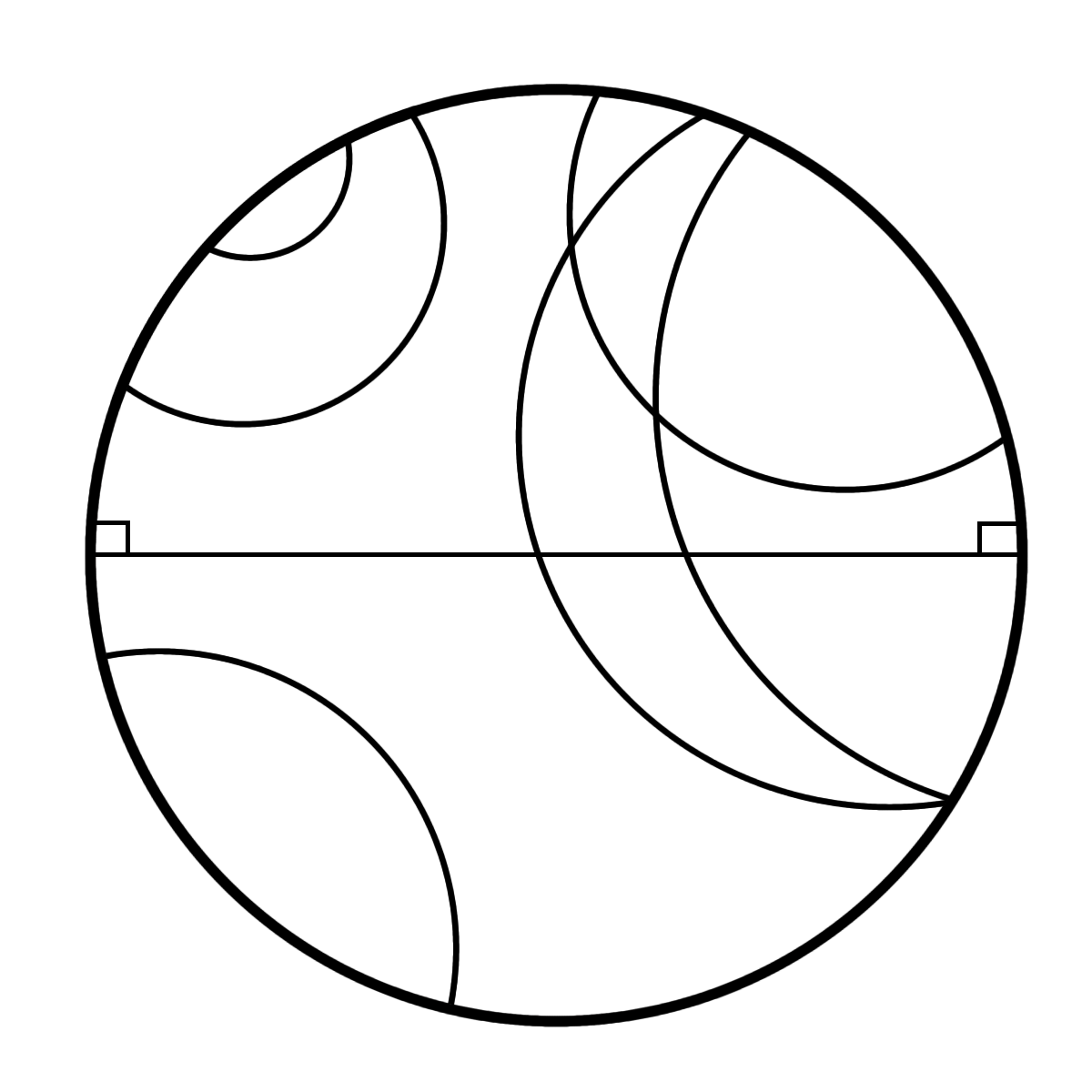}}
    \hspace{2cm}
    {\includegraphics[width=0.35\textwidth]{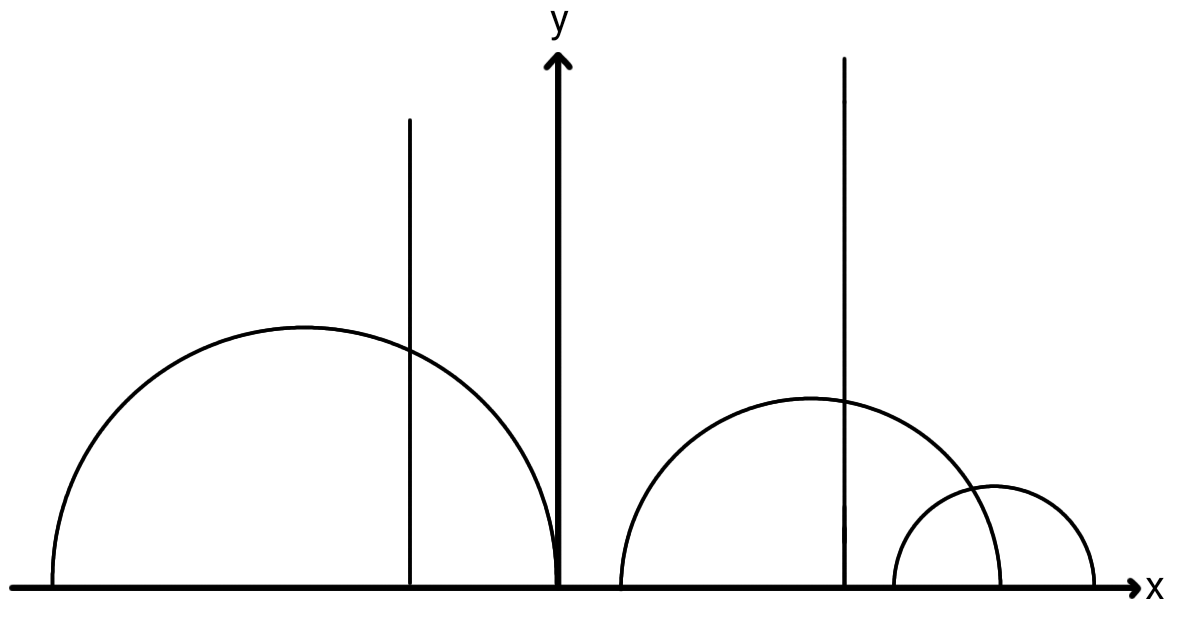}}
    \caption{Geodesics on Poincar\'e disk (left) and semiplane (right).}
\end{figure}

In various models facilitating the comprehension of the Lobachevskian surface, such as the Poincar\'e disk and semiplane, geodesics manifest as curved lines, introducing a degree of counterintuitiveness. The Beltrami-Klein disk model addresses this by rendering geodesics as Euclidean lines, although angle measurement can be intricate. 

In the realm of 3-dimensional models, consider the surface represented by $x^2+y^2-z^2=-k^2$, a hyperboloid of two sheets, and consider only one sheet, named the Weiertrass model. In a pseudo-Euclidean space, a geodesic is the intersection of this surface with a plane passing through the origin  \cite[pp. 188-205]{G}. Another model is the hemisphere model, featuring geodesics in the form of semicircles orthogonal to the equator of the sphere. 

The pseudosphere, introduced by Eugenio Beltrami, is obtained by rotating a tractrix around its asymptote, exhibiting a negative constant Gaussian curvature. Represented as two infinite horns attached at the wide parts, the pseudosphere allows the extension of geodesics in both directions. However, this modification results in a non-simply connected surface, and, according to David Hilbert, in 1901, it was proven that there exists no complete regular surface of constant negative curvature immersed in the three-dimensional Euclidean space, that is, the complete the Lobachevskian surface cannot be represented by a smooth surface with a constant curvature as proposed by Beltrami.

Another tangible 3D model is the crochet model, introduced by Daina Taimina \cite{T}. This woven surface reveals escalating undulations as one moves away from the center. Taimina has proven various properties of this model, enhancing its manipulability, particularly beneficial in an introductory course on non-Euclidean geometry. 
In \cite{O}, the author makes relevant about Taimina's work that the important part is how crocheted model can elucidate mathematical ideas. 

\section{Crocheting a Lobachevskian Surface}
The surface starts with a \emph{magic circle}, which is an adjustable starting round used for crochet patterns that work in crochet rounds, containing a fix number of chains, represented in the middle of Figure \ref{coral} by 6 small circles. The subsequent rows follow a specific pattern, described as follows:
\begin{enumerate}
    \item We begin with $14$ chains to create a manipulable crochet, however, this number can be adjusted as needed. 

\item Over each one of these chains, we knit 3 large chains until completing a row, except in the last one where we knit 4 large chains, obtaining $43$ chains. 
\item For the second row, we continue with the pattern $3332$ over the first $40$ chains, obtaining $110$ chains; and for the remaining $3$ chains, we knit the pattern $3$, for a total of $119$.

\item In the third row, we continue with a pattern $3332$ over the first $116$ chains, obtaining $319$ chains; and for the remaining $3$ chain we knit the pattern $2$, for a total of $325$.
\end{enumerate}

This pattern is designed to create a Lobachevskian surface, taking into account the length of a circle $l = 2\pi \sinh{r}$ as function of its radius. If the magic circle represents $r=1$. In Table \ref{t1} we can see the quantity of chains in each row:
\begin{table}[h]
\centering
\begin{tabular}{| c | c | c |}
\hline
$r$ & $l$ & quantity of chains\\ \hline
$1$ & $7.38$ & $14$\\
$2$ & $22.78$ & $43$\\
$3$ & $62.94$ & $119$\\
$4$ & $171.46$ & $325$\\ \hline
\end{tabular}
\caption{Chains for the knitted surface.}
\label{t1}
\end{table}

\begin{figure}[h]
    \centering
    \includegraphics[width=0.22\textwidth]{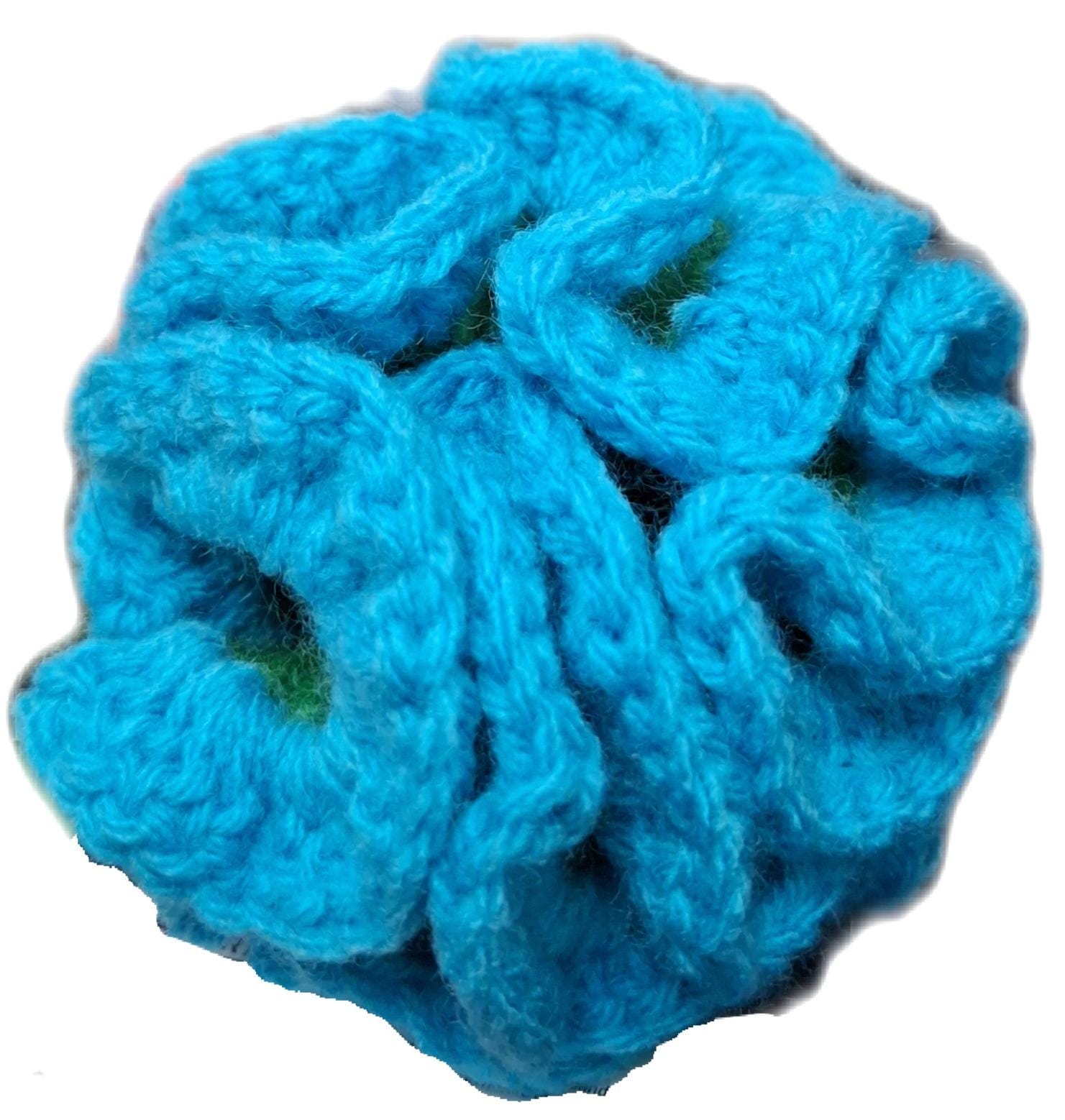}
    \caption{Knitted coral for $r=3$.}
    \label{cc}
\end{figure}

Following this pattern, we obtain the crochet surface in Figure \ref{cc}, which by construction is a Lobachevskian surface.

\section{Parametric Approximation}
Let $v\in [0,2\pi],u\in[0,2]$ and $n$ be a natural number.
The first row ($r=2$) of the knitted coral forms a hyperbolic paraboloid, which has a known parametrization
\begin{align*}
        \textbf{r}  =& u \cos{v} \;\textbf{e}_1 
        + u \sin{v} \;\textbf{e}_2 -u^2 \cos{2v}\;\textbf{e}_3. \end{align*}

When cut, it yields what we denote as a \emph{lettuce}, with parametrization (see Figure \ref{graf1})
\begin{align*}
        \textbf{r}  =& v \; \textbf{e}_1 + u \; \textbf{e}_2-u^2\cos{nv} \; \textbf{e}_3. \end{align*}

\begin{figure}[h]\label{graf1}
    \centering
    \includegraphics[width=0.37\linewidth]{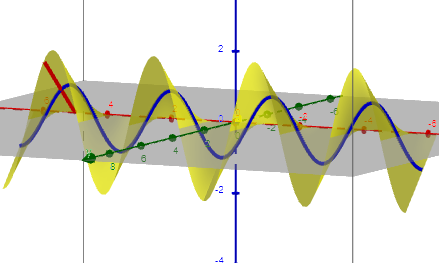}
    \hspace{2cm}
    \includegraphics[width=0.17\linewidth]{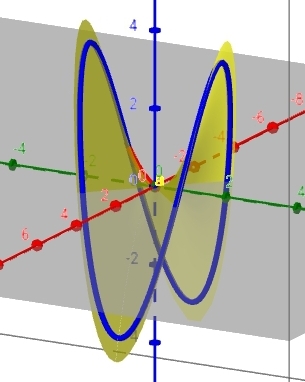}
    \caption{ Geogebra graphics of lettuce (left) and hyperbolic paraboloid (right).}
    \label{graf1}
\end{figure}

\begin{figure}[h]
    \centering
    \includegraphics[width=0.25\linewidth]{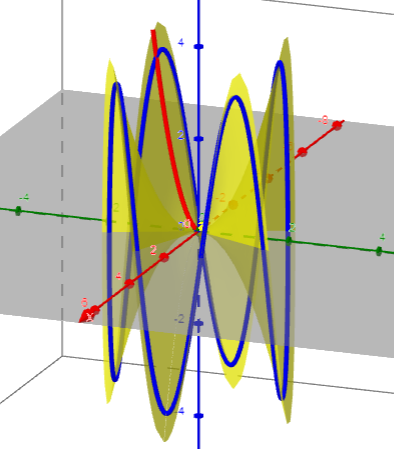}
    \caption{Geogebra graphic of coral.}
    \label{4-coral}
\end{figure}

Drawing upon these two parametrizations, we propose a local parametrization for the knitted coral for $n=4$
\begin{align}\label{coral}
        \textbf{r}  =& u \cos{v} \;\textbf{e}_1 + u \sin{v} \;\textbf{e}_2-u^2 \cos{nv}\;\textbf{e}_3.\end{align}
It is worth mentioning that the $4$-coral graphic serves as a good approximation of the crocheted coral (see Figure \ref{coral} and  \ref{4-coral}). Now, our next step is to calculate the curvature of the $n$-coral surface.

\subsection{Curvature for $n$-coral} 
First, we calculate the  first and second fundamental form of the surface $S$ with local parametrization \eqref{coral} 
with $u\in [0,2]$, $v\in[0,2\pi]$, and $n\in\mathbb N$, $n>1$ a fixed value. Consider the basis
$\left\{\textbf{r}_u,\textbf{r}_v\right\}$ of $T_P(S)$ where
    \begin{align*}
       \textbf{r}_u&=\cos{v}\textbf{e}_1+ \sin{v}\textbf{e}_2 -2u\cos{nv}\textbf{e}_3\\
         \textbf{r}_v&=-u\sin{v}\textbf{e}_1 +u\cos{v}\textbf{e}_2 +nu^2\sin{nv}\textbf{e}_3.
    \end{align*}
The matrix of the first fundamental form of $S$ with respect to the previous basis is$$I=\begin{pmatrix}
        \textbf{r}_u \cdot \textbf{r}_u & \textbf{r}_u \cdot \textbf{r}_v\\
        \textbf{r}_u \cdot \textbf{r}_v & \textbf{r}_v \cdot \textbf{r}_v
    \end{pmatrix}$$
     with
    \begin{align*}
        \textbf{r}_u \cdot \textbf{r}_u = & 1+4u^2\cos^2{nv} \\
        \textbf{r}_u \cdot \textbf{r}_v =& -2nu^3\cos{nv}\sin{nv}\\
        \textbf{r}_v \cdot \textbf{r}_v =& u^2(1+n^2u^2\sin^2{nv}).
    \end{align*}
    
Since
    \begin{align*}
        \textbf{r}_u \times \textbf{r}_v =&  (u^2n\sin{v}\sin{nv}+2u^2\cos{v}\cos{nv}) \;\textbf{e}_1\\
        &+(2u^2\sin{v}\cos{nv} - u^2n \cos{v} \sin{nv})\;\textbf{e}_2\\
        &+u\;\textbf{e}_3\\
        \left\|\textbf{r}_u \times \textbf{r}_v\right\| =& u[n^2u^2\sin^2{nv}+4u^2\cos^2{nv}+1]^{1/2}\\
        =&uA
    \end{align*}
the unit vector field is 
    \begin{align*}
   \textbf{n}=\frac{1}{A}\footnotesize\begin{pmatrix}
            nu\sin{v}\sin{nv}+2u\cos{v}\cos{nv}\\
            2u\sin{v}\cos{nv} - nu \cos{v} \sin{nv}\\
            1
        \end{pmatrix}.
    \end{align*}
    \normalsize
Moreover
    \begin{align*}
        \textbf{r}_{uu}  &= -2\cos{nv}\;\textbf{e}_3\\
        \textbf{r}_{uv}&=-\sin{v}\;\textbf{e}_1+\cos{v}\;\textbf{e}_2+2nu\sin{nv}\;\textbf{e}_3\\
        \textbf{r}_{vv} &= -u\cos{v}\textbf\;{e}_1-u\sin{v}\;\textbf{e}_2+n^2u^2\cos{nv}\;\textbf{e}_3\\
        \textbf{r}_{uu}\cdot \textbf{n} &= -\frac{2}{A}\cos{nv},
        \textbf{r}_{uv}\cdot \textbf{n}  = \frac{nu}{A}\sin{nv},
        \textbf{r}_{vv}\cdot \textbf{n} = \frac{(n^2-2)u^2}{A}\cos{nv}.
    \end{align*}
The matrix  of the second fundamental form of $S$ is
    $$II=\frac{1}{A}\begin{pmatrix}
         -2\cos{nv} & nu\sin{nv}\\
        nu\sin{nv} & (n^2-2)u^2\cos{nv}
    \end{pmatrix}.$$
Since
    \begin{align*}
        I^{-1}&=\frac{1}{\textrm{det}(I)}\begin{pmatrix}
            \textbf{r}_v \cdot \textbf{r}_v & -(\textbf{r}_u \cdot \textbf{r}_v)\\
            -(\textbf{r}_u \cdot \textbf{r}_v) & \textbf{r}_u \cdot \textbf{r}_u\\
        \end{pmatrix}\\
        &=\frac{1}{(\textbf{r}_u \cdot \textbf{r}_u)(\textbf{r}_v \cdot \textbf{r}_v)-(\textbf{r}_u \cdot \textbf{r}_v)^2}\begin{pmatrix}
            \textbf{r}_v \cdot \textbf{r}_v & -(\textbf{r}_u \cdot \textbf{r}_v)\\
            -(\textbf{r}_u \cdot \textbf{r}_v) & \textbf{r}_u \cdot \textbf{r}_u\\
        \end{pmatrix}
    \end{align*}
    the Matrix of the Weingarten transformation is
    \begin{align*}
        &W_p=I^{-1}II\\
        &=\frac{1}{\textrm{det}(I)}\footnotesize{\begin{pmatrix}
         (\textbf{r}_{uu} \cdot \textbf{n})(\textbf{r}_v \cdot \textbf{r}_v)-(\textbf{r}_{uv} \cdot \textbf{n})(\textbf{r}_u \cdot \textbf{r}_v) & -(\textbf{r}_{uu} \cdot \textbf{n})(\textbf{r}_u \cdot \textbf{r}_v)+(\textbf{r}_{uv} \cdot \textbf{n})(\textbf{r}_u \cdot \textbf{r}_u)\\
        (\textbf{r}_{uv} \cdot \textbf{n})(\textbf{r}_v \cdot \textbf{r}_v)-(\textbf{r}_{vv} \cdot \textbf{n})(\textbf{r}_u \cdot \textbf{r}_v) & -(\textbf{r}_{uv} \cdot \textbf{n})(\textbf{r}_u \cdot \textbf{r}_v) +(\textbf{r}_{vv} \cdot \textbf{n})(\textbf{r}_u \cdot \textbf{r}_u)
    \end{pmatrix}}.
    \end{align*}

As a result we obtain the Gaussian curvature  of a $n$-coral, which is
    $$K=-\frac{2(n^2-2)\cos^2{nv}+n^2\sin^2{nv}}{(n^2u^2\sin^2{nv}+4u^2\cos^2{nv}+1)^{3/2}}\cdot$$    

For the $2$-coral the curvature is $\frac{-4}{(4u^2+1)^{3/2}}$, which is constant over circunferences, that is, where $u$ is fixed. As $u$ increases, the curvature approaches zero.
For the $4$-coral is $\frac{12\sin^2{4v}-28}{(12u^2\sin^2{4v}+4u^2+1)^{3/2}}\cdot$ In Table \ref{curvature} we provide some values. 

\begin{table}[h]
\centering
    \begin{tabular}{|c|c|c|c|c|c|c|c|c|}
            \hline
            $v$ & \multicolumn{4}{ c| }{$2\pi$}& \multicolumn{4}{ c| }{$\pi/2$} \\ \hline
            $u$	&  0.5	&  1 &   1.5 & 2  &  0.5	&  1 &   1.5 &   2\\
            \hline
            $\kappa$	& -9.89	& -2.50 &-0.88	& -0.39 & -9.89	& -2.50 &-0.88	& -0.39\\
            \hline
        \end{tabular}
        \caption{Values of the curvature for some values of $4$-coral.}
\label{curvature}
\end{table}

As we obtain a non-constant negative curvature for the $4$-coral, we only achieve an approximation of a Lobachevskian surface through the parametrization of the knitted coral.

\section{Conclusions}
In the literature, a key aspect of the Lobachevskian surface is the concept of curvature, which naively measures how much we deviate from being flat. Various articles, primarily authored by Daina Taimina and her collaborators, explore knitted surfaces to enhance understanding, particularly of the Lobachevskian surface. However, a gap exists regarding the calculation of curvature for these proposed crocheted surfaces.

To achieve this objective, we propose a parametrization of the knitted coral. The resulting object exhibits negative curvature, although it does not remain constant throughout, thus it cannot be deemed a faithful model. Nonetheless, it can be considered an approximation as it captures the deformation obtained through the crochet technique.

\bigbreak

First and second author's Address: km 5 via Puerto Colombia, Barranquilla,\\
Department of Mathematics and Statistic, Universidad del Norte, 080002,
\\
Colombia.\\
First author email: estradaisabella@uninorte.edu.co\\
Second author email: mejiala@uninorte.edu.co\\[4pt]

\label{last}
\end{document}